\documentclass[preprint,12pt]{elsarticle}

\usepackage{xcolor}
\usepackage{graphicx}
\usepackage[version=4]{mhchem}
\usepackage{amssymb}
\usepackage{amsmath}
\def\eps{\ensuremath\varepsilon}
\usepackage{lineno}
\usepackage[margin=1.0in]{geometry}
 \usepackage[colorlinks=true]{hyperref}
 \usepackage{centernot}
 
 \usepackage{amssymb}

\newtheorem{remark}{Remark}
\journal{arXiv}

\usepackage[utf8]{inputenc}
\usepackage{csquotes}

\begin{document}

\begin{frontmatter}

\title{On the validity of the stochastic quasi-steady-state approximation 
in open enzyme catalyzed reactions: \textit{Timescale separation or singular 
perturbation?}}

\author[label1]{Justin Eilertsen}\fnref{cor0}
\author[label1,label2]{Santiago Schnell\fnref{cor1}}

\address[label1]{Department of Molecular \& Integrative Physiology, University of Michigan 
Medical School, Ann Arbor, MI 48109, USA}
\address[label2]{Department of Computational Medicine \& Bioinformatics, University of 
Michigan Medical School, Ann Arbor, MI 48109, USA}

\fntext[cor0]{New Address: Mathematical Reviews, American Mathematical Society, 
416 $4th$ Street, Ann Arbor, MI 48103}
\fntext[cor1]{Corresponding author. New affiliation: Department of Biological Sciences, 
and Department of Applied and Computational Mathematics and Statistics, University of 
Notre Dame, Notre Dame, IN 46556, USA. E-mail: santiago.schnell@nd.edu.}

\begin{abstract}
The quasi-steady-state approximation is widely used to develop simplified deterministic 
or stochastic models of enzyme catalyzed reactions. In deterministic models, the 
quasi-steady-state approximation can be mathematically justified from singular perturbation 
theory. For several closed enzymatic reactions, the homologous extension of the 
quasi-steady-state approximation to the stochastic regime, known as the \textit{stochastic} 
quasi-steady-state approximation, has been shown to be accurate under the analogous 
conditions that permit the quasi-steady-state reduction of the deterministic counterpart. 
However, it was recently demonstrated that the extension of the stochastic quasi-steady-state
approximation to an \textit{open} Michaelis--Menten reaction mechanism is only valid under 
a condition that is far more restrictive than the qualifier that ensures the validity of 
its corresponding deterministic quasi-steady-state approximation. In this paper, we suggest 
a possible explanation for this discrepancy from the lens of geometric singular perturbation 
theory. In so doing, we illustrate a misconception in the application of the 
quasi-steady-state approximation: timescale separation does not imply singular perturbation.
\end{abstract}

\begin{keyword}
Singular perturbation \sep stochastic process, quasi-steady-state approximation, 
Michaelis--Menten reaction mechanism, Langevin equation, linear noise approximation, 
slow scale linear noise approximation, Fenichel theory, Tikhonov's theorem
\end{keyword}

\end{frontmatter}

\newpage

\section{Introduction}\label{sec:0}
The quasi-steady-state approximation (QSSA) is a simple but powerful model reduction 
technique in the field of mathematical biochemistry. In a reduced model of an enzyme 
catalyzed reaction, the QSSA decreases the dimensionality of the full model in both 
state and parameter space which, consequently, increases the tractability of quantifying 
enzyme activity from experimental data or the modeling of the enzyme catalyzed 
reactions \cite{Stroberg2016,Choi2017}. Furthermore, as a course-grained model, the QSSA eases the 
computational burden of numerical simulations by omitting the impact of negligible changes 
in state that take place over extremely small timescales \cite{Sanft,Agarwal2012,Burrage2008}. 

Depending on the specific application, the QSSA comes in two varieties: deterministic 
and stochastic. Regardless of whether or not the deterministic or stochastic approximation 
is employed, the challenge always resides in justifying its use in an application. 
Historically, the deterministic QSSA has been derived from the careful application of 
singular perturbation analysis \cite{Heineken1967,Segel1988,Segel1989}, 
also known as Tikhonov theory \cite{Tikhonov1952} in Eastern Europe, or Fenichel 
theory \cite{Fenichel1979} in the Western hemisphere. The existing link between the QSSA 
and singular perturbation theory is favorable: since singular perturbation theory has 
been developed quite extensively over the past several 
decades \cite{kuehn2015,Wechselberger2020}, there is a rich mathematical foundation 
that supports the legitimization of the QSSA.

In contrast to its deterministic cousin, the justification of stochastic QSSA is 
significantly more difficult, especially in regimes far from the thermodynamic limit 
where the chemical master equation (CME) prevails as the physically relevant model. 
The stochastic multiscale method has provided a rigorous justification of the 
stochastic QSSA~\cite{Ball2006,kang2013,Kang2017}; other justifications have leveraged 
scaling methodologies and singular perturbation methods~\cite{Rao2003,Mastny2007}. Regardless 
of how the stochastic QSSA is justified, there is a developed interest in deriving 
techniques capable of ascertaining, a priori, when the stochastic QSSA will generate 
accurate statistical moments.

Presently, there are several methodologies that rely on the linear noise 
approximation (LNA) that are very good at determining the reliability of the stochastic 
QSSA \cite{Agarwal2012,KIM2014,KIM2015,Thomas2011,Thomas2012,ThomasPO}. In the particular 
case of Michaelis--Menten (MM) type reactions, the stochastic QSSA has been shown to be 
very accurate when the reaction is closed, meaning that it is void of any influx 
(or outflux) of the pertinent chemical species. Perhaps most surprisingly, the conditions 
that support the validity of the stochastic QSSA have been shown to be the 
same\footnote{The use of ``same" here is slightly abusive; deterministic and stochastic 
rate constants differ in terms of their units, and the state space is discrete--as opposed 
to continuous--in the realm of the CME.} conditions 
for the validity of deterministic QSSA~\cite{Sanft}. On the other hand, for open MM-type 
reactions, the proposed qualifier that establishes the accuracy of the stochastic QSSA 
is far more restrictive than the criterion that certifies the effectiveness of the 
corresponding deterministic QSSA \cite{Thomas2011,Othmer2020}. The natural question that 
follows is: \textit{Why}?

The answer may at least partially hinge on the mathematical theory that 
ratifies the legitimacy of the deterministic QSSA. Although it is commonly assumed that 
the reliability of the deterministic QSSA is due to singular perturbation theory, recent 
analyses have clarified that this is not always the 
case \cite{Noethen2009,Noethen2011,Goeke2017}. In fact, in the instance of the open MM 
reaction mechanism, \citet{OpenMMin} demonstrate that there are in 
fact \textit{two} mathematical mechanisms that legalize the QSSA: near invariance of 
the QSS manifold, and singular perturbation reduction. Hence, the validity of the 
deterministic QSSA in the open MM reaction mechanism is not necessarily attributable to 
a singular perturbation. However, when the MM reaction mechanism is closed, the 
QSSA \textit{is} the a result of a singular perturbation. Thus, adding influx to the 
MM reaction mechanism can blur the boundary between singular and non-singular perturbation 
scenarios. This observation is a possible small step in understanding why 
the stochastic QSSA works so well for the closed MM reaction mechanism, but is apparently 
only favorable under more restrictive conditions when open. 

In what follows, we re-examine the role of competing QSSA mechanisms and the possible 
impact this has on the validity of the stochastic QSSA for the open MM reaction mechanism. 
Specifically, we show that the conditions that guarantee the accuracy of the stochastic 
QSSA for the open MM reaction mechanism in the linear noise regime are in agreement with 
a singular perturbation scenario in the corresponding deterministic limit. 

\section{The deterministic formation of the open Michaelis--Menten reaction mechanism}\label{sec:2}
In this paper, we analyze the MM reaction mechanism supplied with a constant influx of 
substrate, which was first examined by \citet{Stoleriu2004}. Let S denote a substrate 
molecule, E an enzyme molecule, C an enzyme-substrate complex, and P a product molecule. 
The reaction mechanism is
\begin{align}\label{mm1}
   \emptyset \ce{->[$k_0$] S},\qquad\ce{S + E <=>[$k_1$][$k_{-1}$] C ->[$k_2$] E + P},
\end{align}
which is known as the open MM reaction mechanism. In contrast, the \textit{closed} 
MM reaction mechanism is void of additional substrate influx, and corresponds to the 
special case when $k_0$ is identically zero.

In the thermodynamic limit, the open MM reaction mechanism (\ref{mm1}) is accurately 
modelled by the following set of rate equations derived from the law of mass action: 
\begin{subequations}\label{MA}
\begin{align}
\dot s&= k_0-k_1(e_T-c)s + k_{-1}c,\label{m1}\\
\dot c&= k_1(e_T-c)s -(k_{-1}+k_2)c\label{m2}.
\end{align}
\end{subequations}
The parameters $k_{1},k_2$ and $k_{-1}$ are deterministic rate constants, $k_0$ is 
the rate of substrate influx, and lowercase $s$, $c$, $p$ denote the respective 
concentrations of substrate, complex, and product. The concentration of E, $e$, can 
be determined a posteriori from the conservation law,
\begin{equation}
e_T = e + c,
\end{equation}
where $e_T$ is the total enzyme concentration, a conserved quantity. Likewise, the 
concentration of P, $p$, can be determined solely from the time course solution of $c$:
\begin{equation}
    \dot{p}=k_2c.
\end{equation}
Since $c\leq e_T$, the limiting rate of the reaction is $v:=k_2e_T$. Furthermore, there 
exists a stable fixed point in the first quadrant, $(s,c):=(\gamma,\nu)$,
\begin{equation}\label{FP}
    \gamma := \cfrac{\alpha K_M}{1-\alpha}, \quad \nu :=\alpha e_T,\quad \alpha:=\cfrac{k_0}{v}
\end{equation}
whenever $k_0 < v$ and the Michaelis constant, $K_M:=(k_{-1}+k_2)/k_1$, is bounded.

To derive the deterministic QSSA from the system~(\ref{MA}), one assumes that after an initial 
transient $\dot{c} \approx 0$, so that the algebraic relationship,
\begin{equation}\label{Qc}
    c=\cfrac{e_Ts}{K_M+s},
\end{equation}
is approximately valid. The curve defined by (\ref{Qc}) is identical to the $c$-nullcline, 
but we will refer to it as the QSS manifold.\footnote{From this point onward, we often 
will refer to the curve defined by (\ref{Qc}) as the QSS manifold.} Insertion of 
(\ref{Qc}) into (\ref{m1}) yields
\begin{equation}\label{sQSSA}
    \dot{s} = k_0 - \cfrac{v\,s}{K_M +s},
\end{equation}
which is known as the \textit{standard} QSSA (sQSSA). 

The relevant question in the application of the sQSSA is always: Under what conditions does 
(\ref{sQSSA}) prevail as a reliable approximation to the full system (\ref{MA})? On the 
deterministic end of the thermodynamic spectrum, the answer to this question for both the open 
and closed MM reaction mechanisms is well-established. As carefully discussed 
in \cite{EILERTSEN2020}, the accuracy of (\ref{sQSSA}) is 
holds as long as\footnote{Locally, near a specific value of $s$, the metric $e_T/(K_M+s)\ll 1$ 
is often sufficient to validate the sQSSA. For more details, please 
consult~\cite{Segel1988,Segel1989,EILERTSEN2020}.}
\begin{equation}\label{segel}
   \eps_{SS} := \cfrac{e_T}{K_M} \ll 1. 
\end{equation}
Colloquially, (\ref{segel}) is known as the \textit{timescale separation} condition. The 
condition (\ref{segel}) was originally justified by Palsson \cite{PALSSON1984,PALSSON1987} via 
linearization of the mass action equations, but the phrase ``timescale separation" follows from 
the work of Segel~\cite{Segel1988} and Segel \& Slemrod~\cite{Segel1989}. Following Segel, but 
using a more rigours scaling, \citet{EILERTSEN2020} demonstrated that 
the sQSSA for the closed MM reaction is valid as long as the timescale that accounts for the 
fast transient, $t_C := 1/(k_{-1}+k_2)$, is much shorter than the timescale 
required for a significant amount of substrate to deplete, $t_S:=1/k_1e_T$:
\begin{equation}
    \eps_{SS} := \cfrac{t_C}{t_S}= \cfrac{e_T}{K_M}.
\end{equation}

\section{The stochastic formulation of the open Michaelis--Menten reaction mechanism}\label{sec:3}
The deterministic mass action equations (\ref{MA}) are only valid in the thermodynamic 
limit, in which both the number of molecules, $n$, and system volume, $\Omega$, approach 
infinity in a way that leaves the respective concentrations of each species finite. Far 
from the thermodynamic limit, where molecular copy numbers are finite, stochastic 
fluctuations persist. In this regime, the master equation, whose solution gives the 
probability of finding the system in state $X$ at time $t$, is the appropriate model to 
employ. The specific CME that corresponds to (\ref{mm1}) is
\begin{multline}\label{CME}
    \cfrac{\partial}{\partial t}P(n_S,n_C,t)=\bigg[\cfrac{k_1}{\Omega}(E_S^{+1}E_C^{-1}n_S(E_T-n_C) + \Omega k_0(E_S^{-1}-1) \\ + k_{-1}(E_S^{-1}E_C^{+1}-1)n_C + k_2(E_C^{-1}-1)n_C\bigg]P(n_S,n_C,t).
\end{multline}
In (\ref{CME}), $P(n_S,n_C,t)$ is the probability of finding the system with $n_S$ substrates 
molecules and $n_C$ complex molecules at time $t$, $\Omega$ is the volume of the system, and 
$E_X^{\pm j}$ are step operators:
\begin{equation}
 E_X^{\pm j} f(X,Y,Z)=f(X\pm j,Y,Z).  
\end{equation}

The corresponding \textit{stochastic} sQSSA is given by
\begin{equation}\label{stochQSSA}
    \cfrac{\partial}{\partial t}P(n_S,t)=\bigg[\Omega k_0 (E_S^{-1}-1) - (E_S^{+1}-1)\cfrac{k_2e_Tn_S}{K_M+n_S/\Omega}\bigg]P(n_S,t).
\end{equation}
A formal derivation of the stochastic QSSA for the closed MM reaction can be found in \cite{Rao2003}; 
the open stochastic QSSA (\ref{stochQSSA}) can be found in \cite{Thomas2011}. In either case, 
the implicit assumption made in the application of the stochastic sQSSA is that the propensity 
function, $a(n_S)$, that gives the probability that a product molecule is formed within the 
infinitesimal window $[t,t+\text{d}t)$ is identically
\begin{equation}\label{prop}
    P(n_S-1, t+\text{d}t|n_S,t):= a(n_S)\text{d}t = \cfrac{k_2e_Tn_S}{K_M+n_S/\Omega}\;\text{d}t.
\end{equation}
Hence, the propensity function for the depletion of substrate is adopted directly from 
the non-elementary rate equation of the deterministic sQSSA. The challenge that emerges 
is thus to determine the specific conditions that permit the use of (\ref{prop}) in 
the Gillespie algorithm. \\

\begin{remark}
In what follows we will be interested in understanding the validity of the stochastic sQSSA near 
the deterministic equilibrium point (\ref{FP}). The reason for this is that, from a computational 
perspective, we can view the open MM reaction (\ref{MA}) as being embedded in a possibly larger 
chemical network. Thus, there is utility in minimizing the computational cost of simulation via 
Gillespie's algorithm by applying the stochastic sQSSA.
\end{remark}

In the linear noise regime, the legitimacy of (\ref{stochQSSA}) in the neighborhood of the
deterministic stationary point (\ref{FP}) was rigorously analyzed by \citet{Thomas2011}. 
These authors discovered that timescale separation is necessary, but not sufficient for 
the validity of (\ref{stochQSSA}) near $(s,c)=(\gamma,\nu)$. Specifically, they demonstrated 
that the variance in the number of substrate molecules under steady--state conditions 
obtained from the full CME (\ref{CME}) can differ from the variance of substrate under 
steady--state conditions computed from the stochastic sQSSA (\ref{stochQSSA}) by as much 
as $30\%$. Given the success of stochastic reductions for the closed reaction~\cite{Sanft,Kang2019}, 
this result is somewhat surprising.

\begin{figure}[htb!]
  \centering
    \includegraphics[width=10.0cm]{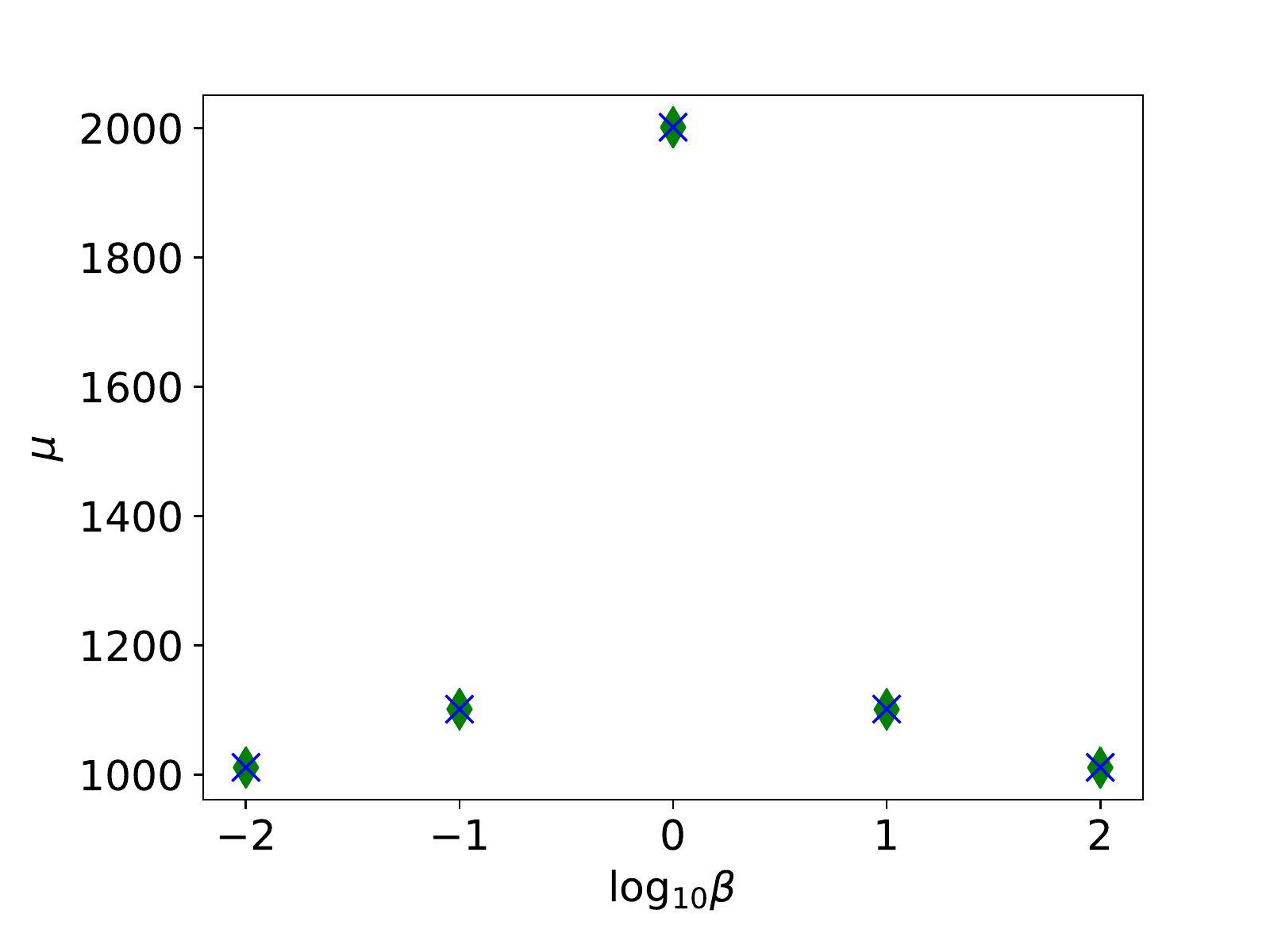}
    \includegraphics[width=10.0cm]{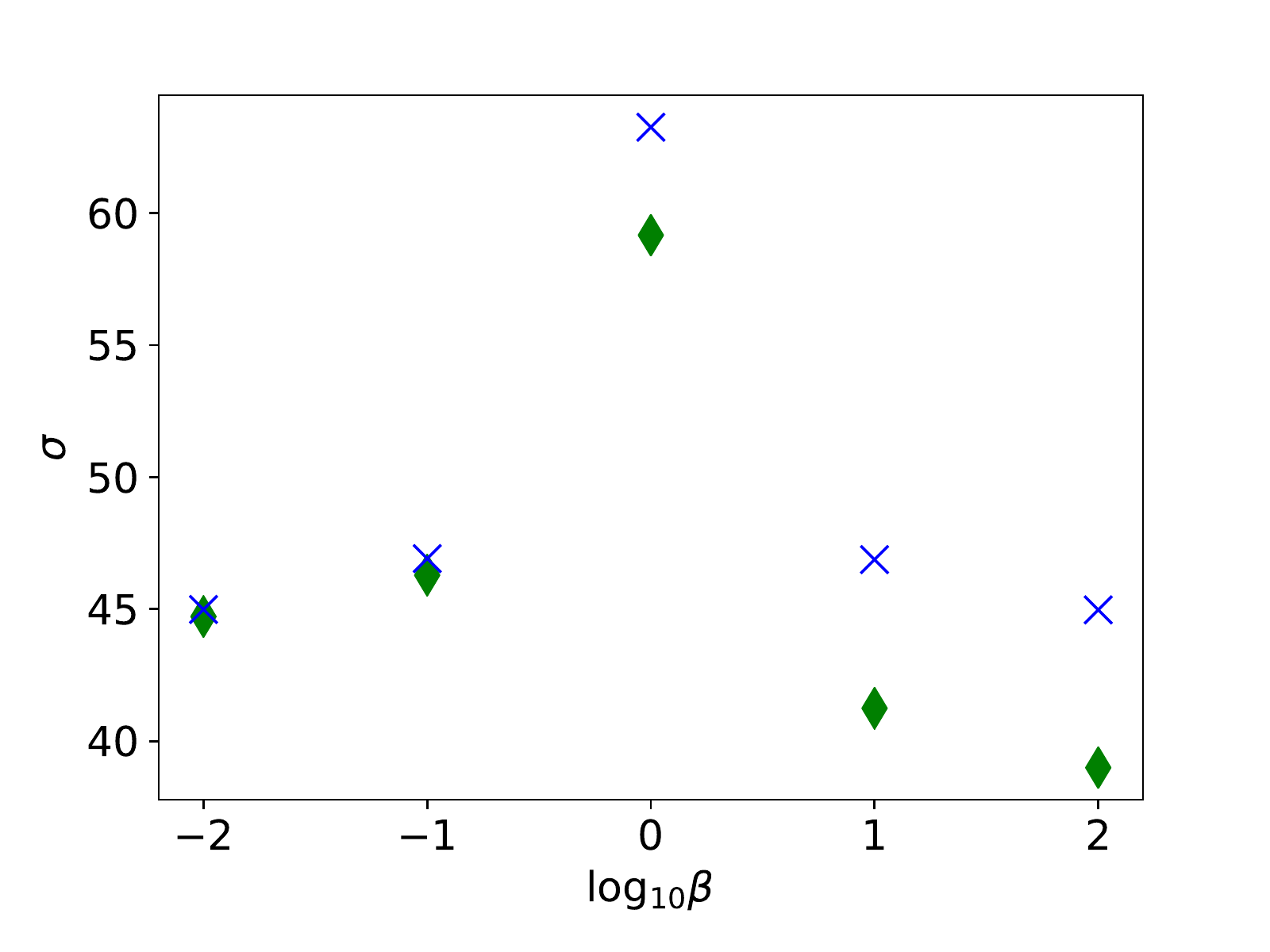}
  \caption{\textbf{Accuracy of the stochastic sQSSA for the open Michaelis--Menten
  reaction mechanism.} If $\alpha \sim 1/2$, then $\beta \ll 1$ is required to ensure the accuracy of (\ref{stochQSSA}). In both panels, $\Omega =1$, the total number of enzyme molecules, $n_E$, is $10$, $k_1=1.0$, $k_0=0.5v$ and $k_{-1},k_2$ are varied so that $k_{-1},k_2 \in [10,100,1000]$ with $\beta =0.01,0.1,1.0,10.0$ and $100.0$. Moreover, in each simulation, parameters are chosen so that $e_T/K_M = 0.01$ is constant. All units are arbitrary. {{\sc top}}: In this panel, $\log_{10} \beta$ is plotted versus the numerically--approximated mean of $n_S$, $\mu$, the number of substrate molecules, under steady--state conditions. The green diamonds are the means estimated from (\ref{CME}) via $10^7$ realizations produced by Gillespie algorithm; the blue crosses are the means estimated from (\ref{stochQSSA}) via $10^7$ realizations generated Gillespie's algorithm. {{\sc bottom}}: In this panel, $\log_{10} \beta$ is plotted versus the numerically--approximated standard deviation of $n_S$, $\sigma$, under steady--state conditions. Green diamonds are standard deviations recovered via realizations of Gillespie's algorithm applied to (\ref{CME}). Blue crosses are standard deviations obtained from realizations generated by Gillespie's algorithm applied to (\ref{stochQSSA}). Clearly, $\beta \ll 1$ is necessary to recover accurate standard deviations when $\alpha = 1/2$.
 } \label{FIG1}
\end{figure}

In order to derive a qualifier that yields a useful a priori indication of the accuracy 
of (\ref{stochQSSA}), Thomas et al. \cite{Thomas2011} analyzed the open MM reaction 
mechanism in the linear noise regime. Performing a $\Omega$-expansion on (\ref{CME}) 
and (\ref{stochQSSA}) and discarding terms that are $o(\Omega^{-1/2})$ yields their 
respective LNAs~\cite{VKX}. From the LNAs, expressions for the substrate variance under
steady--state conditions can computed by solving a Lyapunov equation. 
From \citet{Thomas2011}, these are:
\begin{subequations}
\begin{align}
 \sigma^2_{\text{full}} &= \cfrac{\gamma}{\Omega}\cdot\bigg(1+\cfrac{\gamma}{K_M}\cdot\cfrac{K_S+\gamma}{K_M+\gamma}\cdot\cfrac{1}{1+\varepsilon}\bigg), \qquad K_S:=k_{-1}/k_1,\label{var1}\\
  \sigma^2_{\text{red.}} &= \cfrac{\gamma}{\Omega}\cdot\bigg(1+\cfrac{\gamma}{K_M}\bigg),\label{var2}
\end{align}
\end{subequations}
where $\varepsilon$ is given by:
$\varepsilon:=(e_T-\nu)/(K_M+\gamma)\leq \varepsilon_{SS}$, $\sigma^2_{\text{full}}$ is the substrate 
variance under steady--state conditions obtained from the LNA to full CME (\ref{CME}), 
and $\sigma^2_{\text{red.}}$ is the substrate variance under steady--state conditions 
obtained from the LNA to (\ref{stochQSSA}). By equating $\varepsilon$ with zero 
in (\ref{var1}), \citet{Thomas2011} found that, for $0 \leq \alpha < 1$, it 
is \textit{also} necessary that
\begin{equation}\label{add}
    \bigg|\cfrac{\sigma_{\text{full}}^2-\sigma_{\text{red.}}^2}{\sigma_{\text{full}}^2}\bigg|=\cfrac{(1-\alpha)\alpha\beta}{1+\beta[1-\alpha(1-\alpha)]} \ll 1, \quad \beta := \cfrac{k_{2}}{k_{-1}}, \quad \alpha := \cfrac{k_0}{k_2e_T}\equiv \cfrac{k_0}{v}
\end{equation}
hold in order to ensure the accuracy of the stochastic sQSSA near the deterministic 
stationary point (\ref{FP}). From inspection of (\ref{add}), we observe that 
$\sigma_{\text{red.}}^2$ is a good approximation to $\sigma_{\text{full}}^2$ when 
$0 \leq \alpha \ll 1$ or as $\alpha \to 1$. However, when $\alpha \sim 1/2$, it is 
necessary that $\beta \ll 1$; see {{\sc figure}} \ref{FIG1}. Note that these restrictions 
are in addition to, and distinct from, the timescale separation condition $\eps_{SS}\ll 1$.
Furthermore, while (\ref{add}) was derived under the assumption that the reaction is 
occurring in a region within the thermodynamic spectrum where LNA is valid (i.e., very 
close to the thermodynamic limit), numerical simulations suggest that, when combined 
with timescale separation, (\ref{add}) is a good a priori indicator of the accuracy 
of (\ref{stochQSSA}); again, see {{\sc figure}} \ref{FIG1}. With that said, let us 
state the following:
\begin{remark}
We are not suggesting that the conditions $\varepsilon_{SS} \ll 1$ and (\ref{add}) are 
sufficient for the validity of (\ref{stochQSSA}), as they were obtained by analyzing 
the LNAs and not the CME. However, based on reports in the literature~\cite{ssLNA}, the 
validity of the reduced LNA does appear to yield a good a priori indication of the 
reduced CME's validity.
\end{remark}

The additional constraint defined by (\ref{add}) lies in stark contrast to validity of the 
stochastic QSSA applied to the \textit{closed} MM reaction mechanism (i.e., when $k_0$ is 
identically zero). Earlier studies demonstrated, with great clarity, that the stochastic sQSSA 
for the closed MM reaction is valid under the same conditions that ensure the accuracy of the 
deterministic sQSSA \cite{Sanft}. Moreover, analyses of the deterministic open MM equations 
indicate that the timescale separation condition (\ref{segel}) is also sufficient for the validity 
of the open sQSSA \cite{Stoleriu2004,OpenMMin}. But, \textit{why} does 
the addition of one elementary reaction (i.e., substrate influx) give rise to a more restrictive 
set of qualifiers for the validity of the stochastic sQSSA? After all, the qualifier 
$e_T/k_M \ll 1$ is independent of the amount of substrate present in the system, and thus 
it seems reasonable to assume that the addition of substrate should not compromise the 
integrity of the stochastic QSSA. Nevertheless, as \citet{Thomas2011} have clearly shown, it 
does. Ultimately, the answer to this question may partially lie in the application of singular 
perturbation theory, the basics of which are reviewed in Section~\ref{sec:spt}.

\section{Singular perturbations, near invariance, and the validity of the stochastic 
quasi-steady-state approximation}\label{sec:4}
In this section we review the basics of Fenichel theory, and demonstrate how QSSAs can be 
recovered from slow manifold projection.

\subsection{Singular perturbation theory: The basics} \label{sec:spt}
A perturbed dynamical system may be expressed in the general form
\begin{equation}\label{dy}
    \dot{z}= w(z) + \eps G(z,\eps)
\end{equation}
 where $0 < \eps \ll 1$ is a parameter that is proportional to the ratio of fast and slow 
 timescales. The type of perturbation, i.e., regular or singular, is determined by the number of 
 stationary points in the unperturbed vector field, $w(z)$. If, in the singular limit that 
 coincides with $\eps =0$, there exists a compact and differentiable manifold, $M$, 
 consisting of non-isolated stationary solutions,
 \begin{equation}
     M \subseteq \mathcal{S}:=\{z\in \mathbb{R}^n: w(z)=0\},
 \end{equation}
then the perturbation is \textit{singular}, and $M$ is referred to as the \textit{critical 
manifold}. Moreover, if the critical manifold is normally hyperbolic (and attracting so that 
the real parts of the non-zero eigenvalues of the Jacobian at $z$, $Dw(z)$, are strictly 
less than zero) such that
\begin{equation}\label{splitting}
   \mathbb{R}^n := T_z M \oplus N_z, \quad T_zM=\ker Dw(z), \;\;\text{with rank} \;\;Dw(z)=n-\dim M,\;\; \forall z \in M,
\end{equation}
where $T_zM$ is the tangent space of $M$ and $N_z$ is complementary\footnote{$N_z$ is 
the range, $\mathcal{R}$, of the Jacobian, $Dw(z)$, for $z\in M$.} to $T_z M$, then 
Fenichel theory~\cite{Fenichel1979} says that the perturbed dynamical system (\ref{dy}) 
is approximated via projection of the perturbation, $G(z,0)$,
onto $T_zM$:
\begin{subequations}
\begin{align}
    \Pi^{M}:\mathbb{R}^n &\mapsto T_zM \quad \forall z \in M,\\
    \dot{z} &= \Pi^M  G(z,0).
    \end{align}
\end{subequations}
We refer the reader to \ref{AppendA} as well as \cite{Goeke2017,Wechselberger2020,OpenMMin} 
for details regarding the explicit construction of $\Pi^M$, which is straightforward 
to compute for the open MM reaction.

\subsection{Tikhonov--Fenichel Parameter Values}
In order to justify the application of any form of the QSSA\footnote{For a review 
describing the different forms of the QSSA, we invite the reader to 
consult \cite{Schnell:2003:CEK,EILERTSEN2020}} from singular perturbation theory, it 
imperative to identify the specific conditions under which the vector field of interest 
possesses a critical manifold. Tikhonov--Fenichel parameter 
values (TFPV) provide the necessary conditions for the existence of a normally hyperbolic 
critical manifold in chemical systems \cite{Goeke2017}. For the open MM reaction mechanism,
let $\pi\in \mathbb{R}^5_+$ denote the parameter vector: 
$\pi :=[k_0 \;\; e_T\;\; k_1 \;\; k_2\;\;k_{-1}]^T$. The interesting TFPVs and their 
corresponding critical manifolds, $M$, are as follows:
\begin{subequations}\label{Cman}
\begin{align}
\pi_1&=[0\;0\;k_1\;k_2\;k_{-1}], \quad \& \quad M_1:=\{(s,c)\in \mathbb{R}^2: c=0\},\label{TFPV1}\\
\pi_2&=[0\;e_T\;0\;k_2\;k_{-1}], \quad \& \quad M_2=M_1:=\{(s,c)\in \mathbb{R}^2: c=0\},\label{TFPV2}\\
\pi_3&=[0\;e_T\;k_1\;0\;k_{-1}], \quad \& \quad M_3:=\{(s,c)\in \mathbb{R}^2: c=k_1e_Ts/(k_{-1}+k_1s)\}.\label{TFPV3}
\end{align}
\end{subequations}

A small perturbation to a TFPV results in the validity of a specific QSSA. For example, the TFPV 
$\pi_1$ corresponds to \textit{both} $e_T$ and $k_0$ vanishing in the singular limit. Requiring 
$e_T$ and $k_0$ to be small by defining $e_T\mapsto \eps e_T^*$ and $k_0 \mapsto \eps k_0^*$ with 
$0< \eps \ll 1$ allows us to recast (\ref{MA}) into the perturbation form,
\begin{equation}
    \begin{pmatrix} \dot{s} \\ \dot{c}\end{pmatrix} = \begin{pmatrix}k_1s +k_{-1} \\
    -k_1s -k_{-1}-k_2\end{pmatrix}c + \eps \begin{pmatrix} k_0^*-k_1e_T^*s\\k_1e_T^*s \end{pmatrix},
\end{equation}
with $w(s,c)$ given by
\begin{equation}\label{pert1}
    w(s,c):=\begin{pmatrix}k_1s +k_{-1} \\
    -k_1s -k_{-1}-k_2\end{pmatrix}c,
\end{equation}
and $G(s,c,\eps)$ give by
\begin{equation}\label{pertG}
\eps G(s,c,\eps):=\eps \begin{pmatrix} k_0^*-k_1e_T^*s\\k_1e_T^*s \end{pmatrix}.
\end{equation}
We will refer to (\ref{pert1}) as the \textit{unperturbed vector field} associated with 
the sQSSA, and to (\ref{pertG}) as its associated perturbation. Projecting the perturbation, 
$G(s,c,0)$, onto the tangent space of the critical manifold $M_1$, $T_zM_1$, yields 
the sQSSA (again, see \ref{AppendA} for details):
\begin{equation}
    \dot{s} = k_0 - \cfrac{v\,s}{s+K_M}.
\end{equation}

In a similar manner, defining \textit{both} $k_1$ and $k_0$ to be small 
($k_1 \mapsto \eps k_1^*$ and $k_0\mapsto \eps k_0^*$), or both $k_2$ and $k_0$ to be 
small ($k_2 \mapsto \eps k_2^*$ and $k_0\mapsto \eps k_0^*$) leads to either the small-$s$ 
linear limit of the sQSSA (\ref{lin})
\begin{equation}\label{lin}
   \dot{s} = k_0 - \cfrac{v\,s}{K_M},
\end{equation}
or the quasi-equilibrium approximation (QEA),
\begin{equation}\label{QEA}
\dot{s}=(k_{-1}+k_1s)\cdot\frac{k_0(k_{-1}+k_1s)-k_2k_1e_Ts}{k_1k_{-1}e_T+(k_{-1}+k_1s)^2}, 
\end{equation}
respectively. Again, see \ref{AppendA} of this manuscript, as well as~\cite{OpenMMin} 
for details.

\subsection{Near invariance: A minimal requirement for the justification of the 
quasi-steady-state approximation}
Singular perturbation theory is not the only way to justify a QSSA. In fact, as first 
shown by \citet{Schauer1979}, the ``near invariance" of the 
corresponding QSS manifold (again, usually the $c$-nullcline) is often sufficient to justify 
a QSSA. For example, consider the resulting set of equations that emerge from (\ref{MA}) 
when $e_T = 0$:
\begin{subequations}\label{MA2}
\begin{align}
\dot s&= k_0+k_1cs + k_{-1}c,\label{m12}\\
\dot c&= -k_1cs -(k_{-1}+k_2)c\label{m22}.
\end{align}
\end{subequations}
In the limiting vector field defined by (\ref{MA2}), the manifold $c=0$ is invariant: 
any trajectory starting on $c=0$ stays on it for all time. Moreover, 
$\dot{c}$ is identically zero on $c=0$. Thus, it is heuristically\footnote{In fact, one 
can rigorously make this claim by analyzing the Lie derivative and demanding that the 
vector field be nearly tangent to the QSS manifold; 
see \cite{Schauer1979,Noethen2011,Goeke2017,OpenMMin}.} reasonable to assume that, as 
long as $e_T$ is sufficiently small, the sQSSA should provide a reasonably good 
approximation to the full system (\ref{MA}), since $\dot{c}$ will be close to zero 
when $e_T$ is nearly zero. \citet{Schauer1979}, \citet{Noethen2009,Goeke2017} and 
\citet{OpenMMin} have 
shown that the $c$-nullcline is \textit{nearly invariant} when $e_T$ is 
sufficiently small, which means that the sQSSA (\ref{sQSSA}) is 
\textit{almost}\footnote{We say ``almost" because a solution would require the QSS 
manifold to be perfectly invariant.} a solution to the full system (\ref{MA}). 

While small $e_T$ is sufficient for the validity of the sQSSA, it is insufficient for 
the application of singular perturbation theory. Taking $e_T \to 0$ does not result 
in the formation of a critical manifold of stationary points.  Hence, one cannot 
attribute the validity of the sQSSA to singular perturbation theory in 
regions where only (\ref{segel}) stands, since the vector field does \textit{not} 
contain a critical manifold in the singular limit (again, see \cite{OpenMMin}, 
Section 4 for details). Hence,
\begin{equation}
    \text{near invariance} \centernot\implies \text{singular perturbation}. 
\end{equation}

\begin{remark}
The observation that small $e_T$ does not necessarily lead to a singular perturbation 
scenario (i.e., a vector field with a critical set) for the open MM reaction mechanism~(\ref{MA}) is 
in stark contrast to that of the closed MM reaction. The closed form of the reaction 
corresponds to $k_0=0$. Hence, the manifold $c=0$ is filled with non-isolated stationary 
point when $e_T=0$. Consequently, small $e_T$ is enough to justify the sQSSA from 
singular perturbation theory when the reaction is closed.
\end{remark}

From the practical point of view, we must ask: \textit{Does the precise mathematical 
justification of the sQSSA's validity really matter?} In other words, is it really 
important if the sQSSA derives its justification from near invariance of the QSS 
manifold or from a singular perturbation reduction? In the deterministic regime, the 
answer appears to be no. However, at this point it is not clear if same answer holds 
in the stochastic realm if one tries to extend the form of the non-elementary rate 
equations (from the deterministic QSS reduction) to propensity functions.

In conclusion of this section, we emphasize that the reliability of the deterministic 
sQSSA (\ref{sQSSA}) hinges upon the relative ``smallness" of a specific parameter. At 
the moment, the term ``small" has only been defined generically in terms of an abstract 
parameter ``$\eps$." However, we will clarify the notion of what it means to be 
\textit{small} in Sections~\ref{sec:5} and~\ref{sec:6}. 

\section{The validity of the stochastic QSSA: intimations from Fenichel theory}\label{sec:5}
From the perspective of singular perturbation theory, the fact that timescale separation 
``$\eps_{SS}\ll 1$" does not completely ensure the validity of the stochastic sQSSA is now 
obvious: taking $e_T$ to be exceptionally small is not sufficient to justify the sQSSA 
from singular perturbation theory; it is also necessary that $k_0$ be small. After all, 
both $e_T$ \textit{and} $k_0$ must vanish in the singular limit in order to recover 
the appropriate unperturbed vector field (\ref{pert1}) and the critical manifold $c=0$. 
Of course, the notion of small only makes sense from a relative perception. Concentrations are 
small in comparison to other other concentrations, and rates are small in comparison 
to \textit{other rates}. The total enzyme concentration, $e_T$, is defined to be small 
whenever $e_T\ll K_M$. Thus, $e_T$ is only ``small" when it is much less than $K_M$. On 
the other hand, $k_0$ \textit{may} be designated as small if, for example, it is much 
less than the limiting rate, $v$. Thus, we have that 
\begin{equation}\label{sing}
    \eps_{SS} \ll 1 \quad \text{and} \quad \alpha \ll 1
\end{equation}
comprise a pair of qualifiers that are at least consistent with a singular perturbation 
scenario, since  (\ref{sing}) defines a region in parameter space that is close to the 
TFPV defined by (\ref{TFPV1}). 

However, the story does not end with small $\alpha$. Close inspection of (\ref{add}) reveals 
that $\alpha \ll 1$ is sufficient, but not necessary for the accuracy of (\ref{stochQSSA}) 
since, as $\alpha$ increases with $\eps_{SS}$ remaining small, the stochastic sQSSA is only 
valid if $\beta \ll 1$. For example, the function $\alpha(1-\alpha)$ is maximal 
when $\alpha=1/2$, and so the accuracy of the stochastic QSSA (\ref{stochQSSA}) requires 
$k_2 \ll k_{-1}$ in this case. This is because, at the stationary point where $c:=\nu =\alpha e_T$, 
the influx rate, $k_0$, is \textit{small} compared to the maximal disassociation 
rate,
\begin{equation}
    k_0 = \alpha k_2 e_T \ll k_{-1} e_T.
\end{equation}

The restriction on the size of $k_0$ when $\alpha \sim 1/2$ is also evident 
from (\ref{add}). Rewriting the numerator of (\ref{add}) as
\begin{equation*}
    (1-\alpha)\alpha\beta = (1-\alpha)\lambda, \qquad \lambda:=\cfrac{k_0}{k_{-1}e_T},
\end{equation*}
it becomes clear that one way to ensure (\ref{add}) is small when $\alpha \sim 1/2$ 
is to demand $\lambda \ll 1.$ Thus, if $\alpha \sim 1/2$, then (\ref{stochQSSA}) will 
hold if $k_0$ is much less than the \textit{maximum} disassociation rate: $k_{-1}e_T$. 
Hence, we can interpret the additional qualifier (\ref{add}) as a restriction on the 
relative size of $k_0$,
\begin{equation}
    k_0 \ll \max \{ k_2e_T,k_{-1}e_T\},
\end{equation}
which, together with timescale separation, suggests that the stochastic 
QSSA (\ref{stochQSSA}) is accurate as long as \textit{both} $k_0$ and $e_T$ are be 
comparatively small. And, a restriction on the relative sizes of $k_0$ and $e_T$ is 
perhaps more consistent with a singular perturbation scenario than a near-invariance 
scenario in the deterministic regime.

Finally, let us point out that there is yet another way to interpret the requirement 
that $\beta\ll 1$ when $\alpha \sim 1/2$. Note that if $k_0$ \textit{and} $k_2$ are 
small then presumably, in parameter space, we are close to the TFPV given by (\ref{TFPV3}), as 
opposed to (\ref{TFPV1}). As a result, the QEA~(\ref{QEA}) is valid, since it is now 
$k_0$ and $k_2$ that are ``small." Carefully rewriting the QEA yields
\begin{subequations}\label{QEA1}
\begin{align}
\dot{s}&=(k_{-1}+k_1s)\cdot\frac{k_0(k_{-1}+k_1s)-k_1k_2e_Ts}{k_1k_{-1}e_T+(k_{-1}+k_1s)^2},\\[10pt]
 &=(k_{-1}+k_1s)^2\cdot\frac{k_0-\cfrac{k_1k_2e_Ts}{(k_{-1}+k_1s)}}{k_1k_{-1}e_T+(k_{-1}+k_1s)^2},\\[10pt]
 &=\frac{k_0-\cfrac{k_1 v s}{(k_{-1}+k_1s)}}{\cfrac{k_1k_{-1}e_T}{(k_{-1}+k_1s)^2}+1}\label{fQEA}.
\end{align}
\end{subequations}
Moreover, if $\eps_{SS}\ll 1$ \textit{and} $\beta \ll 1$, then
\begin{equation}
    k_1e_T\ll k_{-1},
\end{equation}
and it follows from (\ref{fQEA}) that the QEA is approximately
\begin{equation}
    \dot{s} = k_0 -\cfrac{v s}{K_S+s}.
\end{equation}
Consequently, we recover a \textit{special case} of the QEA when $\eps_{SS}\ll 1$ and 
$\beta \ll 1$, which is equivalent to the sQSSA when $k_2 \ll k_{-1}$. A parameter 
range in which $\alpha \sim 1/2$ and $\beta \sim 1$ would then be the equivalent of 
``no man's land", since we are presumably somewhere between the TFPVs (\ref{TFPV1}) and (\ref{TFPV3}).

The last component of (\ref{add}) concerns the case in which $\alpha \to 1$. However, this is 
a special case, and justification concerning the validity of the (\ref{stochQSSA}) in this 
region is relatively intuitive. Since the steady-state substrate concentration tends to 
infinity as $\alpha \to 1$, the enzyme molecules are, practically speaking, perpetually 
bound to substrate molecules due to the overwhelming disparity between the populations of 
their respective copy numbers and concentrations. Direct approximation of the invariant 
slow manifold and its proximity to the (non-invariant) QSS manifold also supports the 
validity of the sQSSA when $\alpha \to 1$; see Section 5 in \cite{OpenMMin}.

In summary, our analysis indicates that, in order to ensure the accuracy of the stochastic 
sQSSA (\ref{stochQSSA}), it is absolutely necessary to operate in a parameter regime 
where both $k_0$ is small and $\eps_{SS} \ll 1$. Again, the requirement that both 
(\ref{add}) \text{and} (\ref{segel}) hold is compatible with a singular perturbation 
scenario in the thermodynamic limit.

\section{It's all relative: The utility of dimensionless parameters}\label{sec:6}
The difficulty that resides in interpreting (\ref{add}) arises from the fact that it is 
necessary to demand that \textit{two} parameters be small. In \cite{OpenMMin}, we 
pointed out that singular perturbation theory is really only applicable along rays in 
parameter space. For example, one can fix $k_0 = \alpha k_2 e_T$. By invoking such a 
constraint, it ensures that \textit{both} $k_0$ and $e_T$ vanish by setting $e_T=0$. 
This also allows one to sort of make the claim that $k_0 \sim \mathcal{O}(e_T)$ and 
vice versa. However, it is important to understand that in practical applications, 
what we really want know is: how close are we, in parameter space, to the associated 
TFPV of interest, which can be very difficult to ascertain. Almost all of the theoretical
analyses of the closed MM reaction mechanism are implicitly trying to answer this question 
with respect to \textit{only} $e_T$. For the closed reaction there are four parameters, 
$\pi \in \mathbb{R}^4 := [e_T,k_1,k_2,k_{-1}]^T$, and the validity of the closed sQSSA 
requires $\pi$ to be sufficiently close to the TFPV $\pi_{\text{stand.}} := [0,k_1,k_2,k_{-1}]^T$. 
Thus, the closed sQSSA requires $e_T$ to be small, but what does that mean? Well, 
the answer to this question required a significant amount of
research~\cite{PALSSON1984,Segel1988,Segel1989,EILERTSEN2020,Heineken1967}, the corpus 
of which indicates that $e_T/k_M \ll 1$ is sufficient for the validity of the closed 
sQSSA. The important takeaway from the literature is that the ``smallness" of $e_T$ 
is measured in a comparative sense: $e_T$ is \textit{small} if it is much less than 
$K_M$, and when this holds we can make the argument $\pi$ is close to $\pi_{\text{stand.}}$. 
Thus, the proximity of $\pi$ to $\pi_{\text{stand.}}$ is defined by the magnitude of 
a dimensionless parameter, and it is this distinguishing feature is \textit{crucial}. 
As an example, consider a reaction in which $e_T$ is small (in an arbitrary unit), 
but with a much smaller $K_M$ (in the same unit). The nearest point in parameter space 
would be $\pi=\pi_{\text{rev.}} :=[e_T,k_1,0,0]$. The corresponding vector field,
\begin{subequations}
\begin{align}
\dot{s} &= -k_1(e_T-c)s,\\
\dot{c} &= \;\;\;k_1(e_T-c)s,
\end{align}
\end{subequations}
is equipped with the critical set, $\mathcal{S}$:
\begin{equation}
 \mathcal{S}:= \mathcal{S}_1\cup \mathcal{S}_2 =\{(s,c)\in \mathbb{R}^2: c=e_T\}\; \cup  \;\{(s,c)\in \mathbb{R}^2: s=0\},
\end{equation}
but $\mathcal{S}$ is notably different than $M_1$. Furthermore, $\mathcal{S}$ is not 
even a manifold. It is possible to define compact submanifolds of $\mathcal{S}$ that do 
not contain the intersection $\mathcal{S}_1 \cap \mathcal{S}_2$, and  in fact 
the \textit{reverse} QSSA of the closed MM reaction,
\begin{equation}\label{rev}
    \begin{pmatrix} \dot{s} \\ \dot{c} \end{pmatrix} = \begin{pmatrix}0\\-k_2c_2\end{pmatrix},
\end{equation}
is recovered by projecting the appropriate perturbation onto 
$T_z\mathcal{M}$,\footnote{The ``$\delta$" is included in the definition on $\mathcal{M}$ 
so that the compactness requirement of Fenichel theory is satisfied.}
\begin{equation*}
   \mathcal{M} \subset \mathcal{S}_2 := \{(s,c)\in \mathbb{R}^2: s=0, 0 \leq c \leq e_T - \delta\}, \quad \text{for some} \quad 0 < \delta  < e_T.
\end{equation*}
However, (\ref{rev}) is very different from the closed sQSSA, even though the dimensional 
magnitude of $e_T$ might be small. Thus, the dimensional magnitude of parameter is not 
necessarily sufficient to determine the validity of a QSS reduction.

The need for a qualifier in the form of a dimensionless parameter is not limited to 
chemical kinetics. The ``$\eps$" of Tikhonov and Fenichel theory is colloquially defined as a 
\textit{dimensionless} timescale ratio. Moreover, mathematical timescale ratios are ultimately 
eigenvalue ratios, and the ratio will be dimensionless, even if the eigenvalues 
obtained from the Jacobian of a physical application carry units. Thus, by basing the size 
of a parameter on a dimensional magnitude, we can can easily 
find ourselves in a situation where: (i) there is not a sufficient gap in the 
eigenspectrum, or (ii) our vector field is not close to the TFPV that corresponds to 
the particular QSSA of interest. Hence the need 
for an appropriate \textit{dimensionless parameter}.\footnote{Fenichel 
theory \cite{Fenichel1979} ensures the \textit{eventual} validity of (\ref{sQSSA}) 
as $e_T$ goes to zero with $k_0$ constrained to lie along the ray $k_0=\alpha e_Tk_2$, 
provided initial conditions are sufficiently close to the 
critical manifold and all other parameters are bounded away from zero. However, neither 
TFPV theory nor Fenichel theory indicate how small $e_T$ must be to obtain an accurate 
QSS reduction. Moreover, Fenichel applies only to \textit{compact}, normally hyperbolic 
critical manifolds.} 

In the case of the open MM reaction mechanism (\ref{mm1}), the question is more complex, 
since there is an additional parameter, $k_0$. Thus, one is tasked with ensuring that 
both $e_T$ and $k_0$ are sufficiently small but, small compared to \textit{what}? In 
the deterministic regime, we do not really have to answer the question pertaining to 
\textit{what $k_0$ must be smaller than} to determine when the sQSSA is valid, and 
this is both a pro and a con. It is a pro since all we need to do is vary one parameter 
to ensure the validity of the sQSSA; the accuracy of the open sQSSA is not necessarily 
penalized if we wander too far away from $\pi_1$ in parameter space. It is a con 
because it is difficult to establish \textit{why}, mathematically, the sQSSA is valid. 
Of course, one can constrain $k_0$ and $e_T$ to lie on an appropriate parametric 
curve, but in applications the reaction unfolds at a single point in parameter space, 
and to really attribute the validity of the sQSSA to either a singular perturbation 
scenario or a near invariance scenario, it seems necessary to have some idea about 
the proximity of $\pi$ to $\pi_1$. As we have suggested, the notion of nearness, as 
it pertains to points in parameter space, is best defined in terms of suitable 
dimensionless ratios.

\section{Discussion}\label{sec:8}
In this paper we reiterated how, in the deterministic regime, the validity of the sQSSA 
of the open MM reaction mechanism is attributable to (i) the near invariance of the QSS
manifold, or (ii) a singular perturbation reduction. To determine the mathematical
justification of the sQSSA, it is necessary to have a metric --- defined in terms of a
dimensionless parameter --- that determines the juxtaposition to the TFPV $\pi_1$. A 
singular perturbation scenario for the open MM reaction mechanism requires both 
$e_T$ and $k_0$, the influx rate, to be small. But, small compared to what? Since 
the validity of the LNA to the stochastic QSSA appears to require 
$k_0 \ll \max\{v, k_{-1}e_T\}$, a condition which places a upper bound on the 
relative size of $k_0$, we speculate that this restriction is consistent with 
a singular perturbation scenario. Thus, singular perturbation scenarios in the 
deterministic realm may facilitate a more favorable environment for the validity 
of QSS reductions in the stochastic regime. Let us stress however, that we are 
by no means suggesting that a singular perturbation scenario automatically leads 
the validity of a stochastic QSSA. Several counterexamples that disprove this 
hypothesis already exist in the literature; see~\cite{Janssen1989}. Furthermore, 
we wish to emphasize that our observation that singular perturbation 
scenarios may create more favorable conditions for the accuracy of the stochastic
QSSA (as opposed to just near invariance of the QSS manifold) is certainly not the 
final answer concerning the general validity of stochastic reductions. However, our 
suggestion that the type of mathematical mechanism that permits the deterministic 
reduction may have ramifications concerning the validity of homologous reductions 
in stochastic regime deserves further investigation. 

\appendix
\section{Fenichel theory: projecting onto a slow manifold} \label{AppendA}
In this brief appendix, we briefly describe how to construct $\Pi^M$, the projection 
operator. For more technical details, we invite the reader to 
consult~\cite{Wechselberger2020,Goeke2017,Goeke2015}. Following~\citet{Wechselberger2020}, 
we begin 
with the perturbation form
\begin{equation*}
 \dot{z}= w(z) + \eps G(z,\eps),\quad z\in \mathbb{R}^n.
 \end{equation*}
and let $M$ be a compact subset of  
$w(z)=0$ that forms a $k$-dimensional manifold ($k<n$) comprised of stationary 
points such that: 
\begin{itemize}
    \item For all $z\in M$, the algebraic and geometric multiplicities of the zero 
    eigenvalues of $Dw(z)$ are equal with $T_zM = \ker Dw(z)$ and $\dim \ker Dw(z)= k$, \\
    \item If $\lambda_j$ is a nonzero eigenvalue of $Dw(z)$, then 
    $\mathfrak{Re}(\lambda_j) < 0\quad \forall z\in M$,
\end{itemize}
then, $M$ is locally attracting and there exists a splitting
\begin{equation}\label{splitting-app}
    \mathbb{R}^n:=T_zM \oplus N_z \quad \forall z \in M,
\end{equation}
where $T_zM$ is the tangent space of $M$ at $z$ given by $\{ x \in 
\mathbb{R}^n: x \in \ker Dw(z)\}$, and $N_z$ is the complement to $T_zM$ and 
coincides with range of the Jacobian, $\mathcal{R}(Dw(z))$. The objective from 
the point forward will be to exploit the splitting~(\ref{splitting-app}) and construct 
an \textit{oblique}\footnote{The projection operator $\Pi^M$ is oblique since the 
$Dw(z)$-invariant subspaces, $N_z=\;\text{range}\; Dw(z), \;\;z\in M, \;$ and 
$\ker Dw(z)=T_zM, \;\;z\in M$, are not necessarily orthogonal to one another. 
 } projection operator:
\begin{equation*}
    \Pi^M: \mathbb{R}^n \mapsto T_zM, \qquad I-\Pi^M: \mathbb{R}^n \mapsto N_z.
\end{equation*}

To construct $\Pi^M$, we invoke the factorization
\begin{equation}\label{factor}
    w(z) = P(z)f(z),
\end{equation}
where $P(z)$ is a rectangular matrix function, and the zero level set of $f(z)$ coincides 
with the critical manifold, $M$. The derivative of $f(z)$, $Df(z)$, has full 
rank $\forall z \in M$ and:
\begin{itemize}
    \item The columns of $P(z)$ span the range of the Jacobian $\mathcal{R}(Dw(z))=N_z$, for all $z\in M$.\\
    \item The rows of $Df(z)$ span the orthogonal complement of $\ker Dw(z)=(T_zM)^{\perp}$ for all $z \in M$.
\end{itemize}
It follows that 
since $T_zM$ and $N_z$ are complementary subspaces, the matrix
\begin{equation}\label{projO}
 \Pi^M := I- P(Df P)^{-1}Df
\end{equation}
defines the oblique projection onto $T_zM$. To leading order in $\eps$, 
the reduced flow on the slow manifold is
\begin{equation}
    \dot{z} = \Pi^M G(z,0),
\end{equation}
which \textit{is} the QSS reduction obtained from Fenichel theory~\cite{Fenichel1979}.

\begin{remark}
The critical manifold is said to be normally hyperbolic if the matrix $DfP$ (the eigenvalues of 
which are the non-trivial eigenvalues of $Dw(z)$) is hyperbolic $\forall z\in M$, meaning that 
$\mathfrak{Re}(\lambda_i) \neq 0.$ In most applications, we are interested in \textit{attracting} 
critical manifolds, and thus we often require $\mathfrak{Re}(\lambda_i) < 0.$ However, normally 
hyperbolic critical manifolds are be repelling if $\mathfrak{Re}(\lambda_i) > 0$, or of saddle 
type if the real parts of the eigenvalues of $DfP$ are both positive and negative.
\end{remark}

As an example, we will explicitly compute the sQSSA for the open MM reaction 
mechanism. In standard form, we have
\begin{equation}
    \begin{pmatrix} \dot{s} \\ \dot{c}\end{pmatrix} = \begin{pmatrix}k_{-1} +k_1s  \\
    -k_{-1}-k_2-k_1s \end{pmatrix}c + \eps \begin{pmatrix} k_0^*-k_1e_T^*s\\k_1e_T^*s \end{pmatrix}
\end{equation}
so that $f(s,c)=c$ and 
\begin{equation}
P \equiv \begin{pmatrix}k_{-1} + k_1s  \\
     -k_{-1}-k_2-k_1s\end{pmatrix}.
\end{equation}
The derivative of $f(s,c)$ is $[0 \;\;1]$, and thus $(DfP)^{-1}$ is
\begin{equation}
    -\cfrac{1}{k_1(K_M+s)},
\end{equation}
which is of course a scalar. The product, $PDf$, is given by
\begin{equation}
    \begin{pmatrix} 0 & k_1(K_S+s)\\ 0 & -k_1(K_M+s)\end{pmatrix}, \quad K_S:=k_{-1}/k_1, \quad K_M:=(k_{-1}+k_2)/k_1.
\end{equation}
Computing $\Pi^M$ from (\ref{projO}) yields
\begin{equation}
    \Pi^M:=\begin{pmatrix} 1 & \cfrac{(K_S+s)}{(K_M+s)}\\ 0 & 0\end{pmatrix}.
\end{equation}
To recover the sQSSA, we simply project the perturbation onto the tangent space 
of the critical manifold:
\begin{equation}
   \begin{pmatrix} \dot{s} \\ \dot{c}\end{pmatrix} = \Bigg[\begin{pmatrix} 1 & \cfrac{(K_S+s)}{(K_M+s)}\\ 0 & 0\end{pmatrix}\begin{pmatrix} k_0^*-k_1e_T^*s\\k_1e_T^*s \end{pmatrix}\Bigg]\Bigg|_{c=0} =\begin{pmatrix} k_0^*-\cfrac{k_1k_2 e_T^*s}{k_{-1}+k_2+k_1s}\\ 0\end{pmatrix}.
\end{equation}
The QEA (\ref{QEA}) that corresponds to small $k_0$ \textit{and} small $k_2$ 
is computed in a similar manner. For more details see \cite{OpenMMin}.

\section*{Acknowledgments}
Justin Eilertsen was partially supported by the University of Michigan Postdoctoral Pediatric 
Endocrinology and Diabetes Training Program ``Developmental Origins of Metabolic 
Disorder'' (NIH/NIDDK grant: T32 DK071212).


\end{document}